\documentclass[11pt]{amsart}
\usepackage{amssymb}
\usepackage{color}

\newtheorem{theorem}{Theorem}[section]

\theoremstyle{definition}

\newtheorem{remark}[theorem]{Remark}

\def\Cay{{\rm Cay}}
\def\Aut{{\rm Aut}}
\def\Sym{{\rm Sym}}

\newcommand{\V}{\mathrm{V}}
\newcommand{\Z}{\mathrm{Z}}
\newcommand{\ZZ}{\mathbb Z}

\newcommand{\vGa}{\vec{\Gamma}}

\title[Arc-transitive  digraphs with blocks of given size]{Arc-transitive  digraphs of given out-valency and with blocks of given size}
\author{Luke Morgan}
\author{Primo\v{z} Poto\v{c}nik}
\author{Gabriel Verret}

\address{Luke Morgan\\
Centre for the Mathematics of Symmetry and Computation, School of Mathematics and Statistics (M019)\\
The University of Western Australia\\
Crawley, 6009\\
Australia} 
\email{luke.morgan@uwa.edu.au}

\address{Primo\v{z} Poto\v{c}nik, Faculty of Mathematics and Physics, University of Ljubljana, Jadranska 21, SI-1000 Ljubljana, Slovenia.\newline
\indent Also affiliated with: Institute of Mathematics, Physics and Mechanics, Jadranska 19, SI-1000 Ljubljana, Slovenia.
}
\email{primoz.potocnik@fmf.uni-lj.si}

\address{Gabriel Verret\\
Department of Mathematics, The University of Auckland, Private Bag 92019, Auckland 1142, New Zealand.}
\email{g.verret@auckland.ac.nz}

\begin{document}
\begin{abstract}
Given integers $k$ and $m$, we construct a $G$-arc-transitive graph of valency $k$ and an $L$-arc-transitive 
oriented digraph of out-valency $k$ such that $G$ and $L$ both admit  blocks of imprimitivity of size $m$.
\end{abstract}
\maketitle
\section{Introduction}

A  digraph $\Gamma=(V,A)$ is a set of \emph{vertices} $V$ and a set of \emph{arcs} $A \subseteq V\times V$.  An \emph{automorphism} of $\Gamma$ is a permutation of $V$ that preserves $A$. 
 The digraph $\Gamma$ is \emph{$G$-arc-transitive} if $G$   is a group of automorphisms of $\Gamma$ that is transitive on  $V$ and on $A$. Observe that, in this case, the binary relation $A$ is either symmetric or asymmetric, and we call $\Gamma$ a \emph{graph} or an \emph{oriented digraph}, respectively. The \emph{out-valency} of  $\Gamma$ is then $|\{v \in V \mid (u,v) \in A\}|$ (which is independent of the choice of $u\in V$), and we simply call  this the  \emph{valency} when $\Gamma$ is a graph.  In this short note, we prove the following:

\begin{theorem}
\label{main}
Let $k\geqslant 2$ and  
$m\geqslant 1$
 be integers. 
\begin{enumerate}
\item There exists a finite connected $L$-arc-transitive oriented digraph of out-valency $k$  such that $L$ admits blocks of imprimitivity of size $m$.
\item There exists a finite connected $G$-arc-transitive graph of valency $k$ such that $G$ admits blocks of imprimitivity of size $m$.
\end{enumerate}
\end{theorem}

This answers a question that was posed to us by R\"{o}gnvaldur G. M\"{o}ller and Sara Zemlji\v{c}. They needed examples of such digraphs for a problem about arc-types of infinite vertex-transitive digraphs \cite{moller}.

In Section~\ref{sec:heis groups} we  construct infinite Cayley digraphs  of prime out-valency which are the parents of all our examples. In Sections~\ref{sec:odd m} and ~\ref{sec:even m} we quotient these infinite examples by an appropriate normal subgroup to obtain \emph{finite} Cayley digraphs with the required properties. We then complete the proof in Section~\ref{sec:proof}.

\section{Some infinite Cayley graphs}
\label{sec2}

\subsection{Discrete Heisenberg groups}
\label{sec:heis groups}
 Given  elements $x$ and $y$  of a group, we write   $[x,y]$ for  $x^{-1}y^{-1}xy$ and $x^y$ for $x^{-1}yx$.

Let $n\geqslant 1$ be an integer and let $k=2n+1$. The \emph{discrete Heisenberg group in dimension $k$}, is 
\begin{align*}
H_{k}=\langle x_1,\ldots,x_{2n},z \mid \> &  [z,x_i]=1 \textrm{ for } i\in \{1, \ldots, 2n\}, \\
&[x_i,x_j]=1 \textrm{ if } |j-i|\neq n, \\
&[x_i,x_{i+n}]=z  \textrm{ for } i\in \{1, \ldots, n\} \rangle.
\end{align*}

Note that $H_{k}$ is  generated by $\{x_1,\ldots,x_{2n}\}$ (but giving a ``name'' to $z$ greatly simplifies the notation). Note also that  $\Z(H_{k})=\langle z\rangle$ is infinite cyclic and $H_{k}/\langle z\rangle$ is free abelian of rank $2n$. In particular, every element of $H_k$ can be written uniquely in the form  
$x_1^{\alpha_1}\cdots x_{2n}^{\alpha_{2n}}z^{\alpha_{2n+1}}$, with $\alpha_1,\ldots,\alpha_{2n+1}\in \ZZ$.

Let $t$ be the unique automorphism of $H_k$ such that 
\begin{align*}
x_i^t&=x_i^{-1} \quad\mathrm{ for } \quad i\in\{1,\ldots,2n\},\\ 
z^t&=z. 
\end{align*}
To see that this indeed defines an automorphism, one must first observe that the images of the generators of $H_k$ under $t$
generate $H_k$, and then check that  every defining relation of $H_k$, after substituting each generator in the relation  by its  image under $t$, becomes  a relation in $H_k$. 
For example, the relation $[x_i,x_{i+n}] = z$ becomes $[x_i^{-1},x_{i+n}^{-1}] = z$. Now, in $H_k$, we see that $[x_i^{-1},x_{i+n}^{-1}] =
 x_{i+n}x_i[x_i,x_{i+n}]x_i^{-1}x_{i+n}^{-1} = x_{i+n}x_i z x_i^{-1}x_{i+n}^{-1} = z$, and so the afore mentioned equation is indeed a relation in $H_k$.
 The other defining relations are easier to check and are left to the reader.

Note that $t^2=1$ and let
$$R_k=H_k\rtimes \langle t\rangle.$$

From now on, we assume that $k=2n+1$ is prime. We define an automorphism $b$ of $H_k$ by the following rule:

$x_i^b=\begin{cases}
x_{i+1}x_{n+1} & \textrm{ if } i\in\{1,\ldots,n-1\},\\
x_{i+1} & \textrm{ if } i\in\{n,\ldots,2n-1\},\\
x_1^{-1}x_{n+1}^{-1}\cdots x_{2n}^{-1}  & \textrm{ if } i=2n.\\
\end{cases}$
In \cite[4.1]{Burgiser}, it was shown that $b$ indeed defines an automorphism  of $H_k$ of order $k$. 
(More specifically, the subgroup of $\mathrm{Aut}(H_k)$ acting trivially on $\mathrm Z (H_k)$ contains a subgroup isomorphic to $\mathrm{Sp}(2n,\ZZ)$, and $b$ is simply the automorphism of $H_k$ induced by a matrix of $\mathrm{Sp}(2n,\ZZ)$ that is shown in loc.~cit.~to have order $k$.) Note that $z^b=z$.
We now extend $b$ to an automorphism of $R_k$ by setting
$$t^b :=x_{n+1}t.$$
We need to check that this indeed defines an automorphism of $R_k$, and so we must check that $b$ preserves the defining relations of $R_k$ (in fact, merely those that involve $t$). Note that $(t^b)^2=(x_{n+1}t)^2=x_{n+1}tx_{n+1}t=x_{n+1}x_{n+1}^{-1}=1$. Moreover, for $i\in\{1,\ldots,n-1\}$, we have 
$$(x_i^b)^{t^b}=(x_i^b)^{x_{n+1}t}=(x_{i+1}x_{n+1})^{x_{n+1}t}=(x_{i+1}x_{n+1})^t=x_{i+1}^{-1}x_{n+1}^{-1}=(x_i^b)^{-1}.$$
For $i\in\{n,\ldots,2n-1\}$, 
$$(x_i^b)^{t^b}=(x_i^b)^{x_{n+1}t}=(x_{i+1})^{x_{n+1}t}=(x_{i+1})^{t}=(x_{i+1})^{-1}=(x_i^b)^{-1},$$
and, finally,
\begin{align*}
(x_{2n}^b)^{t^b}&=(x_{2n}^b)^{x_{n+1}t}=(x_1^{-1}x_{n+1}^{-1}\cdots x_{2n}^{-1})^{x_{n+1}t}\\
&=(x_1^{-1}x_{n+1}^{-1}\cdots x_{2n}^{-1}z^{-1})^{t}=x_1x_{n+1}\cdots x_{2n}z^{-1}\\
&=x_{2n}\cdots x_{n+1}x_1=(x_1^{-1}x_{n+1}^{-1}\cdots x_{2n}^{-1})^{-1}\\
&=(x_{2n}^b)^{-1}.
\end{align*}
Thus $b$ extends to an automorphism of $R_k$ (which we also denote by $b$), as claimed. By induction on $i$, one can easily show that, for every $i\in \{1,\ldots, n\}$, we have
$$
t^{b^i} = x_{n+i}\ldots x_{n+1} t \quad \hbox{ and } \quad t^{b^{n+i}} = x_i^{-1}t.
$$
Recall that $b$ induces an automorphism of order $k$ on $H_k$ and, by the above, we see that
 $t^{b^k} = (x_n^{-1}t)^b = (x_n^{-1})^bt^b = t$. Hence $b$, as an automorphism of $R_k$, also has order $k$. Let 
$$ G_k=R_k\rtimes \langle b\rangle,$$
$$L_k = \langle H_k,b\rangle= H_k \rtimes \langle b \rangle,$$
$$S_{k}=\{ t,t^b,\ldots,t^{b^{2n}}\},$$ 
$$P_k = \{x_{n+1}, (x_{n+1})^b,\ldots, (x_{n+1})^{b^{2n}} \}.$$ 
Since the order of $b$, as an automorphism of $R_k$, is $k$, we see that both $P_k$ and $S_{k}$ are closed under conjugation by $\langle b\rangle$. Since $(x_{n+1})^b\neq x_{n+1}$, $t^b\neq t$ and $b$ has prime order, we have
$$|P_k|=|S_{k}|=|b|=k.$$

We claim that $\langle S_{k}\rangle =R_k$. First, note that $t^bt=x_{n+1}\in \langle S_k\rangle$.  For $i\in \{-1,0,\ldots,n-1\}$, we have $x_{n+1}^{b^i}=x_{n+1+i}$ and thus, since $\langle S_k \rangle$ is closed under conjugation by $b$, it follows that $x_n,\ldots,x_{2n}\in \langle S_k\rangle$. Note that $(x_{n-1})^b=x_{n}x_{n+1}\in \langle S_k\rangle$ and thus $x_{n-1}\in\langle S_k\rangle$. Repeating this argument, we get that $x_1,\ldots,x_{n-1}\in\langle S_k\rangle$ and thus  $\langle S_k\rangle =R_k$. This calculation also shows that $\langle P_k \rangle = H_k$.

Let $\Gamma_k=\Cay(R_k,S_k)$ (that is, $\Gamma_k$ is the digraph with vertex set $R_k$ and arc set $\{ (g ,sg) \mid g\in R_k, s \in S_k\}$).
Since $\langle S_k\rangle = R_k$, we see that $\Gamma_k$ is connected and, since
the elements of $S_k$ are involutions, $\Gamma_k$ is a graph of valency $k$.

Let $\vGa_k = \Cay(H_k,P_k)$. Since $\langle P_k \rangle = H_k$, $\vGa_k$ is connected.
Observe that no element of $P_k$ is the inverse of another element in $P_k$, implying that
$\vGa$ is an oriented digraph of out-valency $k$. 

Since $\langle b \rangle$ is a group of automorphisms of $R_k$ that preserves and acts transitively on $S_k$, the group $G_k$, being equal to $R_k\rtimes \langle b \rangle$, acts arc-transitively as a group of automorphisms on the Cayley graph $\Gamma_k$ (see, for example, \cite{godsil}). Similarly, since $\langle b\rangle$ is a subgroup of $\Aut(H_k)$ acting transitively on $P_k$, $L_k$ acts arc-transitively on the digraph $\vGa_k$.

\subsection{$m$ odd}
\label{sec:odd m} 
In this section, we assume that $m \geqslant 3$ is some odd positive integer. Let $N_k=\langle x_1^m,\ldots,x_{2n}^m,z^m\rangle\leqslant H_k$. Note that
$$(x_i^m)^{x_{i+n}}=(x_i^{x_{i+n}})^m=(x_iz)^m=x_i^mz^m\in N_k$$
and, similarly, 
$(x_{i+n}^m)^{x_i}  = (x_{i+n} z^{-1})^m = x_{i+n}^m z^{-m}   \in N_k$ and thus $N_k$ is a normal subgroup of $H_k$. Clearly, $(N_k)^t=N_k$. Moreover, for $i\in \{1,\ldots,2n-1\}$, we see that $(x_i^m)^b$ is equal to either $x_{i+1}^mx_{i+n}^m$ or $x_{i+1}^m$, and is thus an element of $N_k$. On the other hand,
\begin{equation}\label{eq:1}
(x_{2n}^m)^b=(x_{2n}^b)^m=(x_1^{-1}x_{n+1}^{-1}\cdots x_{2n}^{-1})^m=x_1^{-m}x_{n+1}^{-m}\cdots x_{2n}^{-m}z^{-\binom{m}{2}}.
\end{equation}
Since $m$ is odd, $m$ divides $\binom{m}{2}$ and thus $(x_{2n}^m)^b\in N_k$.
 It follows that $N_k$ is a normal subgroup of $G_k$. 
 
Since each $s\in S_k$ is an involution, the element $sN_k \in R_k/N_k$ is also an involution. Hence $\Gamma_{k,m}:=\Cay(R_k/N_k, \{sN_k \mid s \in S_k \})$ is a connected $G_k/N_k$-arc-transitive graph of valency $k$. 

Note that $(x_{n+1})^2 \notin  N_k$, and hence $x_{n+1}N_k \in H_k/N_k$ is not an involution. It follows that
$ \vGa_{k,m}:=\Cay(H_k/N_k, \{sN_k \mid s \in P_k \})$ is a connected $L_k/N_k$-arc-transitive  oriented digraph of out-valency $k$.

Note that $\langle z\rangle N_k / N_k \cong \langle z \rangle/ (\langle z \rangle \cap N_k) = \langle z \rangle / \langle z^m \rangle$ has order $m$, acts semiregularly on the vertices of  $\Gamma_{k,m}$ as well as 
of $\vGa_{k,m}$ and is normal in $G_k/N_k$. Thus the orbits of $\langle z\rangle N_k / N_k$ on the vertices of $\vGa_{k,m}$, respectively, $\Gamma_{k,m}$ are blocks of size $m$ for $L_k/N_k$, respectively, $G_k/N_k$.

Finally, note that $\langle z\rangle N_k / N_k$ is a central cyclic subgroup of $H_k/N_k$ of order $m$, and  that the quotient $(H_k/N_k)/(\langle z\rangle N_k / N_k)$ is isomorphic to $\ZZ_m^{2n}$. It follows that $|\V(\vGa_{k,m})|=|H_k/N_k|=m^k$ and $|\V(\Gamma_{k,m})|=2m^k$.

\subsection{$m$ even}
\label{sec:even m}
Now, we assume $m\geqslant 2$ is even. Let 
\begin{eqnarray*}
v_n=v_{n+1}= & \cdots & =v_{2n}=1, \\
(v_{n-1},v_{n-2},v_{n-3},v_{n-4},\ldots) & = & (0,1,0,1,\ldots), \\
E_k & = & \langle  x_1^mz^{v_1m/2},\ldots, x_{2n}^mz^{v_{2n}m/2},z^m\rangle;
\end{eqnarray*}
In particular, $v_1 = 1$ if $n$ is even, and $v_1=0$ if $n$ is odd. Again, it easy to check that $E_k$ is a normal subgroup of $H_k$; for example,
$(x_i^mz^{v_im/2})^{x_{i+n}}=(x_i^{x_{i+n}})^mz^{v_im/2}=(x_iz)^mz^{v_im/2}=x_i^mz^{v_im/2}z^m\in E_k$. Note that $(x_i^mz^{v_im/2})^t=x_i^{-m}z^{v_im/2}=(x_i^mz^{v_im/2})^{-1}\in E_k$ and thus $(E_k)^t=E_k$. 

We now check that $(E_k)^b=E_k$.  For $i\in\{1,\ldots,n-1\}$, we have
\begin{align*}
(x_{i}^mz^{v_{i}m/2})^b&=(x_i^b)^mz^{v_{i}m/2}=(x_{i+1}x_{n+1})^mz^{v_{i}m/2}=x_{i+1}^mx_{n+1}^mz^{v_{i}m/2}\\
&=(x_{i+1}^mz^{v_{i+1}m/2})(x_{n+1}^mz^{v_{n+1}m/2})z^{(-v_{i+1}-v_{n+1}+v_{i})m/2}.
\end{align*}
Note that $v_{n+1}=1$ and $v_{i}-v_{i+1}\in \{-1,1\}$, hence $z^{(-v_{i+1}-v_{n+1}+v_{i})m/2}\in \langle z^m\rangle$ and $(x_{i}^mz^{v_{i}m/2})^b\in E_k$. It is easy to check that this also holds for $i\in \{n,\ldots,2n-1\}$. It remains to check the case $i=2n$. By \eqref{eq:1}, $(x_{2n}^m)^b=x_1^{-m}x_{n+1}^{-m}\cdots x_{2n}^{-m}z^{-\binom{m}{2}}$ and thus
\begin{align*}
(x_{2n}^mz^{v_{2n}m/2})^b    &   =   x_1^{-m}x_{n+1}^{-m}\cdots x_{2n}^{-m}z^{-\binom{m}{2}}z^{v_{2n}m/2}\\
&=(x_1^mz^{v_1m/2})^{-1}(x_{n+1}^mz^{v_{n+1}m/2})^{-1}\cdots \\ 
&  \cdots (x_{2n}^mz^{v_{2n}m/2})^{-1}z^{-\binom{m}{2}+(v_1+v_{n+1}+\cdots +v_{2n}+v_{2n})m/2}.
\end{align*}
Observe that $v_1+v_{n+1}+\cdots+v_{2n}$ is even. Finally, we have $\binom{m}{2}\equiv m/2 \pmod m$ hence $(x_{2n}^mz^{v_{2n}m/2})^b\in E_k$. This completes our proof that $E_k$ is a normal subgroup of $G_k$. 

As in the previous section, we conclude that $\Cay(R_k/E_k,\{s E_k \mid s\in S_k\})$ is a connected $G_k/E_k$-arc-transitive graph of valency $k$  and order $2m^k$,  and  the orbits of $\langle z \rangle E_k / E_k$ form blocks of size $m$ for $G_k/E_k$.

Recalling the definition of $E_k$, we can show that $x_{n+1}E_k \in H_k/E_k$ has order $2m \geqslant 4$.  It follows that $\Cay(H_k/E_k,\{s E_k \mid s\in P_k\})$ is a connected $L_k/E_k$-arc-transitive oriented digraph of out-valency $k$  and order $m^k$, and  the orbits of $\langle z \rangle E_k / E_k$ form blocks of size $m$ for $L_k/E_k$.

\section{Proof of Theorem~\ref{main}}
\label{sec:proof}
Let $\Gamma=(V,A)$ be a connected $G$-arc-transitive digraph of out-valency $b$. The Cartesian $a$th power of $\Gamma$, denoted $\Gamma^a$, has vertex set $V^a$ and arc set $\{((v_1,\dots,v_k),(u_1,\dots,u_k)) \mid (v_i,u_i) \in A$ for some $i$ and $v_j=u_j$ for all $j\neq i\}$. The digraph $\Gamma^a$ is a connected  
$(G \wr \Sym(a))$-arc-transitive digraph of out-valency $ab$. Moreover, if $G$ admits a block of size $m$, say $B$, then $(G \wr \Sym(a))$ also does, namely the ``diagonal'' block $\{(u,\ldots,u)\mid u\in B\}$. Note that $\Gamma^a$ is an oriented digraph, respectively, graph, if and only if $\Gamma$ is an oriented digraph, respectively, graph.  This reduces the proof of Theorem~\ref{main} to the case of prime out-valency.

We now prove Theorem~\ref{main} in this case. Let $m \geqslant 2$ (the case that $m=1$ is vacuous) and let $k$ be a prime. If $k$ is odd (so $k=2n+1$ for some integer $n\geqslant 1$) the required digraph is shown to exist in Section~\ref{sec:odd m} or ~\ref{sec:even m} depending on the parity of $m$. Thus we may suppose that $k=2$. A cycle of length $2m$ is an arc-transitive graph of valency two with blocks of size $m$.   For $m\geqslant 3$ let $A_m=\ZZ_3\times\ZZ_m$ and  $\vGa_{2,m}=\Cay(A_m,\{(1,1),(1,-1)\})$. There is an automorphism $\tau$ of $A_m$ swapping $(1,1)$ and $(1,-1)$. It follows that $\vGa_{2,m}$ is  an $\langle A_m,\tau\rangle$-arc-transitive  oriented digraph of valency two and order $3m$ admitting blocks of size $m$  (the orbits of a $\ZZ_m$ subgroup of $A_m$). 
For $m=2$, $\vGa_{2,2}=\Cay(\ZZ_4\times\ZZ_2,\{(1,0),(1,1)\})$ is a $\mathrm D_8 \times \ZZ_2$-arc-transitive oriented digraph of valency two, order 8 and with blocks of size 2.

 \section{Concluding remarks} 
\begin{remark}
As remarked earlier, $H_k/N_k$ in Section~\ref{sec:odd m} has centre of order $m$ and the quotient by the centre is isomorphic to $\ZZ_m^{2n}$. In particular, when $m$ is prime, $H_k/N_k$ is an extraspecial group of order $m^k$.
\end{remark}

\begin{remark}
Whilst we have made no effort to compute the (full) automorphism group of the digraphs constructed in this paper, computational evidence suggests that the automorphism group of the graph $\Gamma_{k,m}$ constructed in  Section~\ref{sec2}   is twice as big as $G_k/E_k$.
\end{remark}

\begin{remark}
Our construction may appear  somewhat complicated given the relatively simple question we are trying to answer. Indeed, for some special values of $k$ and $m$, there are   easier constructions. Nevertheless, we have reason to believe that a general solution cannot be much simpler. For example,   consulting Marston Conder's  census of cubic arc-transitive graphs on up to 10,000 vertices \cite{conder}, one can check that   our examples of graphs for $k=3$ in Theorem~\ref{main}(2) are  ``optimal'' (in the sense of having minimal order)  for  $m\in\{5,11,17\}$ and, we believe, for infinitely many other values of $m$.
\end{remark}

\section*{Acknowledgements.}
The authors acknowledge the financial support by the Australian Research Council grant DE160100081 and by
the Slovenian Research Agency (research core funding P1-0294).
 We would also like to thank the Centre for the Mathematics of Symmetry and Computation at The University of Western Australia for its support, as this project started as a problem for the Centre's annual retreat in 2016.   The first and third author thank the second for hosting them at the University of Ljubljana, where much progress on this problem was made. Finally, the first author thanks Marston Conder and Eamonn O'Brien for supporting a visit to The University of Auckland where this work was completed.

\end{document}